\documentclass[12pt,psfig, reqno]{amsart}
\usepackage{amsmath, amsthm, amsfonts, amssymb, latexsym, graphicx}
\usepackage[latin1]{inputenc}

\newtheorem{theorem}{Theorem}[section]

\theoremstyle{definition}
\newtheorem{definition}[theorem]{Definition}
\newtheorem{example}[theorem]{Example}

\numberwithin{equation}{section}
\newcommand{\C}{\mathbb{C}}

 \allowdisplaybreaks
\begin{document}

\title{Fatou and Julia like sets}
\author[K. S. Charak]{Kuldeep Singh Charak}
\address{
\begin{tabular}{lll}
& Kuldeep Singh Charak\\
& Department of Mathematics\\
& University of Jammu\\
& Jammu-180 006\\
& India
\end{tabular}}
\email{kscharak7@rediffmail.com }

\author[A. Singh]{Anil Singh}
\address{
\begin{tabular}{lll}
& Anil Singh\\
& Department of Mathematics\\
& University of Jammu\\
& Jammu-180 006\\
& India
\end{tabular}}
\email{anilmanhasfeb90@gmail.com }

\author[M. Kumar]{Manish Kumar}
\address{
\begin{tabular}{lll}
& Manish Kumar\\
& Department of Mathematics\\
& University of Jammu\\
& Jammu-180 006\\ 
& India\\
\end{tabular}}
\email{manishbarmaan@gmail.com}

\begin{abstract}  For a family of holomorphic functions on an arbitrary domain, we introduce Fatou and Julia like sets, and establish some of their interesting properties.
\end{abstract}

\renewcommand{\thefootnote}{\fnsymbol{footnote}}

\maketitle

\section{Introduction and Main results}
 Throughout, we shall denote by $\mathcal{H}\left(D\right)$ the class of all holomorphic functions on a domain $D\subseteq \mathbb{C}$.
 A subfamily $\mathcal{F}$ of $\mathcal{H}\left(D\right)$ is said to be normal if every sequence in $\mathcal{F}$ contains a subsequence that converges locally uniformly on $D$. $\mathcal{F} $ is said to be normal at a point $z_0\in D$ if it is normal in some neighborhood of $z_0$ in $D$ (see \cite{schiff, zalcman1}).

\medskip

Let $f$ be an entire function and let $f^{n}:= \underbrace{f\circ f\circ\cdots \circ f }_{n-times}(n\geq1)$ be the $n$-th iterate of $f$. The \textit{Fatou set} of $f$, denoted by $F(f)$, is defined as 
$$F(f)=\left\{z\in \mathbb{C}: \left\{f^{n}\right\} \mbox{ is a normal family in some neighborhood of } z \right\}$$ 
and, the complement $\mathbb{C}\setminus F(f)$ of $F(f)$ is called the \textit{Julia set }of $f$ and is denoted by $J(f)$. $F(f)$ is an open subset of $\mathbb{C}$ and $J(f)$ is a closed subset of $\mathbb{C}$, and both are completely invariant sets under $f$. The study of Fatou and Julia sets of holomorphic functions is a subject matter of Complex Dynamics for which one can refer to \cite{Bergweiler-5, gamelin,  steinmetz}.

\medskip

 For a given domain $D$ and a subfamily $\mathcal{F}$ of $\mathcal{H}\left(D\right)$, we denote by $F(\mathcal{F})$, a subset of $D$ on which $\mathcal{F}$ is normal and $J(\mathcal{F}):= D\setminus F(\mathcal{F})$. If $\mathcal{F}$ happens to be a family of iterates of an entire function $f$, then $F(\mathcal{F})$ and $J(\mathcal{F})$ reduce to the Fatou set of $f$ and the Julia set of $f$ respectively, therefore, it is reasonable to call $F(\mathcal{F})$ and $J(\mathcal{F})$ as \textit{Fatou and Julia like sets}. Note that Julia set of an entire function is always non-empty (see \cite{Bergweiler-5}) whereas Julia like set $J(\mathcal{F})$ can be empty. For example, consider the family
 $$\mathcal{F}:=\{f(az+b): a,b \in \mathbb{C}, a\neq 0\},$$
 where $f$ is a normal function on $\mathbb{C}$ (see \cite{schiff}, p. 179). Then since $f$ is a normal function on $\mathbb{C}$, $\mathcal{F}$ is a normal family on $\mathbb{C}$. That is, $F\left(\mathcal{F}\right)=\mathbb{C}$ and hence $J(\mathcal{F})=\phi$. 
	
	Also, it is interesting to note that Julia set of any meromorphic function is an uncountable set (see \cite{Bergweiler-5})  but Julia like set is not so, for example, $J(\mathcal{F})=\{0\}$, where $\mathcal{F}:=\{nz: n\in \mathbb{N}\}\subset \mathcal{H}\left(\mathbb{D}\right)$, where $\mathbb{D}$ is the open unit disk.
	
	\medskip

	 If $\mathcal{F}$ and $\mathcal{G}$ are two subfamilies of $\mathcal{H}\left(D\right)$, then $J(\mathcal{F}\cap \mathcal{G}) \subset J(\mathcal{F})\cap J(\mathcal{G})$, however $J(\mathcal{F}\cap \mathcal{G}) = J(\mathcal{F})\cap J(\mathcal{G})$ may not hold in general. For example, let $$\mathcal{F}=\{nz: n \in \mathbb{N}\}\cup \{n(z-1): n \in \mathbb{N}\}$$ and $$\mathcal{G}=\{n(z-1):n \in \mathbb{N}\}\cup \{e^{nz}: n \in \mathbb{N}\}.$$ be the families of entire functions.	Then $J(\mathcal{F}\cap \mathcal{G})= \{1\}$ and $J(\mathcal{F})\cap J(\mathcal{G})=\{0,1\}$.
 
\medskip

This article is devoted to the problem of normality of families of mappings that have been actively studied recently, (see \cite{mc, krss, mrv, rm, rss}). In particular, we give some interesting properties of Fatou and Julia like sets.

%%%%%%%%%%%%%%%%%%%%%%%%%%%%%%%%%%%%%%%%%%%%%%%%%%%%%%%%%%%%%%%%%%%%%%%%%%%%%%%%%%%%%%%%%%%%%%%%%%%%%%%%%%%%%%%%%%%%%%%%%%%%%%%%%%%%%%%%%%%%%
\begin{theorem}\label{union} (a) If $\mathcal{F}_{1}$ and $\mathcal{F}_{2}$ are two subfamilies of $\mathcal{H}\left(D\right)$, then $J(\mathcal{F}_{1}\cup \mathcal{F}_{2})= J(\mathcal{F}_{1})\cup J(\mathcal{F}_{2})$.
 
(b) If $z_{0}\in J(\mathcal{F})$ and $N$ is any neighborhood of $z_0$, then $\C\setminus U $ contains at most one point, where $U=\bigcup _{f\in \mathcal{F}} f(N)$.
\end{theorem}
%%%%%%%%%%%%%%%%%%%%%%%%%%%%%%%%%%%%%%%%%%%%%%%%%%%%%%%%%%%%%%%%%%%%%%%%%%%%%%%%%%%%%%%%%%%%%%%%%%%%%%%%%%%%%%%%%%%%%%%%%%%%%%%%%%%%%%%%%%%%
\begin{example} For $\alpha\in\mathbb{C}$, consider one-parameter family of entire functions $\mathcal{F}_{\alpha}:=\{n(z-\alpha):n\in\mathbb{N}\}$. Then $\mathcal{F}_{\alpha}$ is not normal at $z = \alpha$, that is, $\mathcal{F}_{\alpha}$ is not normal in any open set containing $z=\alpha$. Consider the family of entire functions $\mathcal{F}=\cup_{|\alpha|\leq 1}\mathcal{F}_{\alpha}$. Then we show that $J(\mathcal{F})=\{z: |z|\leq 1 \}$ and hence $Int \left(J\left(\mathcal{F}\right)\right) \neq \phi$ and $J\left(\mathcal{F}\right)\neq \mathbb{C}$.
 The inclusion $\{z: |z|\leq 1 \}\subset J(\mathcal{F})$ holds trivially. To show the other way inclusion, let $z_0 \in\mathbb{C}$  such that $|z_0|>1$ and let $\{f_n\}$ be a sequence in $\mathcal{F}$. Then we have two cases:\\
{\it Case-I}: When $\{f_n\}$ has a subsequence $\{f_{n_k}\}$ which is locally bounded at $z_0$. 
  
	 In this case by Montel's theorem $\{f_{n_k}\}$ further has a  subsequence which converges uniformly in some neighborhood of $z_0$. Thus, $\{f_n\}$ has a subsequence which converges uniformly in some neighborhood of $z_0$. That is, $z_0\in F(\mathcal{F})$.  \\
{\it Case-II}: When $\{f_n\}$ has no subsequence which is locally bounded at $z_0$. 

   Since $f_n(z)=m_n(z-\alpha_n)$, where $m_n \in \mathbb{N}$ and $|\alpha_n|\leq 1$ for each $n\in\mathbb{N}$,  it follows that $\{m_n\}$ has an increasing subsequence $\{m_{n_k}\}$ which converges to $\infty$. Let $N \subset \{z:|z|>1\}$ be a small neighborhood of $z_0$ . Then $\{f_{n_k}\}$ converges uniformly to $\infty$ in $N$. Thus $\{f_n\}$ has a subsequence which converges uniformly in some neighborhood of $z_0$. That is, $z_0 \in F\left(\mathcal{F}\right)$.
	 
	 Thus in both the cases we find that $J(\mathcal{F})\subset \{z:|z|\leq 1\}$. Hence $J(\mathcal{F})= \{z:|z|\leq 1\}.$
\end{example}
 Note that for a family of iterates of an entire function $f$,  $J(f)=\C$ or $J(f) $ has empty interior \cite[Lemma 3]{Bergweiler-5}.

\medskip

A set $A\subset D$ is said to be  \textit{forward invariant (backward invariant)} under the family $\mathcal{F}$ if, for each $f\in \mathcal{F},$ $f(A)\subset A \left(f^{-1}(A)\subset A\right).$ 

   If $\mathcal{F}_0$ is a semi-group of entire functions, then $F(\mathcal{F}_0)$ is forward invariant and $J(\mathcal{F}_0)$ is backward invariant under the family $\mathcal{F}_0$, (see \cite{Hinkkanen}), whereas for an arbitrary subfamily $\mathcal{G}$ of $\mathcal{H}\left(D\right)$, $F(\mathcal{G})\mbox{ and }J(\mathcal{G})$ may not be forward invariant or backward invariant. For example, $J\left(\mathcal{G}\right)$  is not forward invariant as well as backward invariant, for $\mathcal{G}=\{nz:n\in\mathbb{N}\}\cup\{z^n:n\in\mathbb{N}\}$.  Forward invariance of  $J(\mathcal{F})\mbox{ and }F(\mathcal{F})$  for the family $\mathcal{F},$ implies the following: 
%%%%%%%%%%%%%%%%%%%%%%%%%%%%%%%%%%%%%%%%%%%%%%%%%%%%%%%%%%%%%%%%%%%%%%%%%%%%%%%%%%%%%%%%%%%%%%%%%%%%%%%%%%%%%%%%%%%%%%%%%%%%%%%%%%%%%%%%%%
\begin{theorem}\label{emptyinterior}
Let $\mathcal{F}$ be a subfamily of  $\mathcal{H}\left(D\right)$. Then  the following holds :
\begin{itemize}
 \item[(a)] If $J(\mathcal{F})$ is forward invariant, then $J(\mathcal{F})=D \mbox{ or }Int(J(\mathcal{F}))=\phi$. In particular, if $\mathbb{C}\setminus D$ contains at least two points, then $Int(J(\mathcal{F}))=\phi.$ 
\item[(b)] If $J(\mathcal{F})$ contains at least two points and $F(\mathcal{F})$ is forward invariant, then $J(\mathcal{F})$ is a perfect set. 
\end{itemize}
\end{theorem}
%%%%%%%%%%%%%%%%%%%%%%%%%%%%%%%%%%%%%%%%%%%%%%%%%%%%%%%%%%%%%%%%%%%%%%%%%%%%%%%%%%%%%%%%%%%%%%%%%%%%%%%%%%%%%%%%%%%%%%%%%%%%%%%%%%%%%%%%%%%
\begin{example}\label{forward} 
	Let $\mathcal{F}_1=\left\{nz:n\in\mathbb{N}\right\}$,  $\mathcal{F}_2=\left\{z^n:n\in\mathbb{N}\right\}$. 
	Then $J(\mathcal{F}_1\cup\mathcal{F}_2)=\left\{z:|z|=1\right\}\cup\left\{0\right\}$. Clearly,  $J\left(\mathcal{F}_1\cup\mathcal{F}_2\right)$  is not perfect and  $F\left(\mathcal{F}_1\cup\mathcal{F}_2\right)$ is not forward invariant. 
\end{example}
 Example \ref{forward} shows that the condition, ``$F(\mathcal{F})$ is forward invariant'' in Theorem \ref{emptyinterior} can not be dropped.
 
\medskip

 Recall that a point $z_0\in D$ is said to be a periodic point of an entire function $f$, of order $k$, if $f^{k}(z_0)=z_0$. In the dynamics of transcendental entire functions, it is well known that Julia set is the closure of the repelling periodic points (see \cite{schwick}). This can be extended for Julia like set too. In this context, we need some basic notations from the Nevanlinna value distribution theory of meromorphic functions (see \cite{Hayman-1}).

 Let $f$ be a meromorphic function on $\mathbb{C}$. The proximity function $m(r,a,f)$ of $f$ and the counting function $N(r,a,f)$ of $a-$points of $f(a\neq\infty) $ are given by 
$$m(r,a,f):=\frac{1}{2\pi}\int_{0}^{2\pi}{\log^{+}\frac{1}{|f(re^{i\phi})-a|}}d\phi.$$
 For $a=\infty$, we write 
$$m(r,f):=\frac{1}{2\pi}\int_{0}^{2\pi}{\log^{+}|f(re^{i\phi})|}d\phi.$$
$$N(r,a,f):=\int_{0}^{r}{\frac{n(t,1/f-a)}{t}}dt$$
and $$N(r,f):=\int_{0}^{r}{\frac{n(t,f)}{t}}dt,$$ where $n(t,1/f-a)$ is the number of $a-$points of $f$ in $|z|\leq t$ and in particular, $n(t,f)$ is the number of poles of $f$ in $|z|\leq t$. The characteristic function of $f$, denoted by $T(r,f)$, is given by
$$T(r,f)=m(r,f)+N(r,f)$$
and it behaves like $\log^{+}{M(r,f)}$ whenever $f$ happens to be an entire function, where $M(r,f)=\max_{|z|=r}|f(z)|$.
Further, we define 
$$\delta(a,f)=1-\limsup_{r\rightarrow\infty} \frac{N(r,a,f)}{T(r,f)}$$ 
and is called the Nevanlinna deficiency of $f$ at $a$, and the truncated defect is given by 
$$\Theta(a,f)=1- \limsup_{r\rightarrow\infty}\frac{\overline{N}(r,a,f)}{T(r,f)}$$
where $\overline{N}(r,a,f)$ is the counting function of $f$ corresponding to the distinct $a-$points of $f$, that is, by ignoring the multiplicities of $a-$points of $f$.

%%%%%%%%%%%%%%%%%%%%%%%%%%%%%%%%%%%%%%%%%%%%%%%%%%%%%%%%%%%%%%%%%%%%%%%%%%%%%%%%%%%%%%%%%%%%%%%%%%%%%%%%%%%%%%%%%%%%%%%%%%%%%%%%%%%%%%%%%%%
\begin{theorem}\label{1} Let $\mathcal{F}$ be a family of transcendental entire functions and  $J(\mathcal{F})$ contains at least three points. Then for any $w_0\in J(\mathcal{F})$ and $ f\in\mathcal{F}$ with $\Theta(w_o,f)<\frac{1}{2},$ there exist a sequence $\left\{w_n\right\}\mbox{ such that }w_n\to w_0$ and a sequence $\left\{f_n\right\}\subset\mathcal{F}$ such that $w_n$ is a repelling fixed point of $fof_n.$  
\end{theorem}
The polynomial analogue of  Theorem \ref{1} also holds as follows:
%%%%%%%%%%%%%%%%%%%%%%%%%%%%%%%%%%%%%%%%%%%%%%%%%%%%%%%%%%%%%%%%%%%%%%%%%%%%%%%%%%%%%%%%%%%%%%%%%%%%%%%%%%%%%%%%%%%%%%%%%%%%%%%%%%%%%%%%%%%%
\begin{theorem}\label{poly}
If $\mathcal{F}$ is a family of non-constant polynomials in which for each $w_0\in J\left(\mathcal{F}\right)$, there is $P_0 \in\mathcal{F}$ 
such that $P_0-w_0$ has at least three distinct simple roots. Then $J\left(\mathcal{F}\right)$ is contained in the closure of repelling fixed points of the polynomials of the form $PoQ$, where $P,Q \in\mathcal{F}$.
\end{theorem}
%%%%%%%%%%%%%%%%%%%%%%%%%%%%%%%%%%%%%%%%%%%%%%%%%%%%%%%%%%%%%%%%%%%%%%%%%%%%%%%%%%%%%%%%%%%%%%%%%%%%%%%%%%%%%%%%%%%%%%%%%%%%%%%%%%%%%%%%%%%%
\begin{definition}\label{orbit1}
For a subfamily  $\mathcal{F}$ of $\mathcal{H}\left(D\right)$ and $z\in\mathbb{C},$ define  
\begin{flalign*}\mathcal{O}^{-}_{\mathcal{F}}(z):= &\left\{w\in D:f(w)=z,\mbox{ for some }f\in\mathcal{F}\right\}\\ 
   = & \bigcup_{f\in\mathcal{F}}f^{-1}\{z\} 
\end{flalign*} and  
$$E\left(\mathcal{F}\right):=\{z\in\mathbb{C}:\mathcal{O}^{-}_{\mathcal{F}}(z)\mbox{ is finite}\}.$$
\end{definition}
For a family $\mathcal{F}$ of non-constant entire functions and $z_0\in\mathbb{C}$ , $ \mathcal{O}^{-}_{\mathcal{F}}(z_0)$ is finite implies that $f^{-1}\{z_0\}$ is finite for each $f\in\mathcal{F}.$ In this case $N(r,z_0,f)=O(1)$ and hence $\delta(z_0,f)=1\mbox{ for all }f\in\mathcal{F}.$ While $\delta(z_0,f)=1\mbox{ for all }f\in\mathcal{F}$ may not always imply that $\mathcal{O}^{-}_{\mathcal{F}}(z_0)$ is finite as shown by the following example:
	\begin{example} Let $\mathcal{F}=\left\{\left(z-n\right)e^z:n\in\mathbb{N}\right\}.$ Then $\mathcal{F}$ is a family of transcendental entire functions and  $\mathcal{O}^{-}_{\mathcal{F}}(0)=\left\{n:n\in\mathbb{N}\right\}$ is infinite and $N(r,0,f)=O(\log(r))\mbox{ as }r\to\infty$ and hence $\delta(0,f)=1\mbox{ for all }f\in\mathcal{F}.$
	\end{example}
%%%%%%%%%%%%%%%%%%%%%%%%%%%%%%%%%%%%%%%%%%%%%%%%%%%%%%%%%%%%%%%%%%%%%%%%%%%%%%%%%%%%%%%%%%%%%%%%%%%%%%%%%%%%%%%%%%%%%%%%%%%%%%%%%%%%%%%%%%
 By an extension of Montel's theorem (\cite{cara} p.203), it follows that if  $\mathcal{O}^{-}_{\mathcal{F}}(z_0)$ is omitted by $\mathcal{F}$ on some deleted neighborhood of some  $w\in J(\mathcal{F})$, then $\mathcal{O}^{-}_{\mathcal{F}}(z_0)$  contains at most one point and hence $z_0\in E\left(\mathcal{F}\right).$

%%%%%%%%%%%%%%%%%%%%%%%%%%%%%%%%%%%%%%%%%%%%%%%%%%%%%%%%%%%%%%%%%%%%%%%%%%%%%%%%%%%%%%%%%%%%%%%%%%%%%%%%%%%%%%%%%%%%%%%%%%%%%%%%%%%%%%%%%
Let $\mathcal{F}$ be a uniformly bounded family of holomorphic functions on a domain $D$. Then by Montel's theorem $J(\mathcal{F})=\phi$. Note that $E(\mathcal{F})$ is an infinite set. Indeed, there exists $M>0$ such that $|f(z)|\leq M$ for all $f\in\mathcal{F}$ and so $\{w:|w|>M\}\subset E(\mathcal{F})$ showing that $E(\mathcal{F})$ is uncountable. Let $\mathcal{F}= \{f\in \mathcal{H}(D): f \mbox{ omits two distinct fixed values $a$ and $b$ on } D\}$. Then by Montel's Theorem, $J(\mathcal{F})=\phi$ and $E(\mathcal{F})=\{a,b \}$. The size of $E(\mathcal{F})$ has a definite relation with $J(\mathcal{F})$. In fact, we have the following result:
%%%%%%%%%%%%%%%%%%%%%%%%%%%%%%%%%%%%%%%%%%%%%%%%%%%%%%%%%%%%%%%%%%%%%%%%%%%%%%%%%%%%%%%%%%%%%%%%%%%%%%%%%%%%%%%%%%%%%%%%%%%%%%%%%%%%%%%%%%%%
\begin{theorem}\label{orbit} Let  $\mathcal{F}$ be a subfamily of $\mathcal{H(D)}$.
\begin{itemize}
\item[(a)] If $E \left(\mathcal{F}\right)\neq \phi$, then for $z\notin E \left(\mathcal{F}\right)$, $J(\mathcal{F})\subseteq\overline{\mathcal{O}^{-}_{\mathcal{F}}(z)}.$ 
\item[(b)]If $J(\mathcal{F})\neq \phi$, then $\# E(\mathcal{F})\leq 1.$
\end{itemize}
\end{theorem}
%%%%%%%%%%%%%%%%%%%%%%%%%%%%%%%%%%%%%%%%%%%%%%%%%%%%%%%%%%%%%%%%%%%%%%%%%%%%%%%%%%%%%%%%%%%%%%%%%%%%%%%%%%%%%%%%%%
Following example shows that $E(\mathcal{F})$ may contain exactly one point:
\begin{example}
Let $\mathcal{F}=\{nz: n \in \mathbb{N}\}$ be the family of entire functions. Then $\mathcal{O}_{\mathcal{F}}^{-1}(0)=\{0\}$, it follows that $0 \in E(\mathcal{F})$. Note that $J(\mathcal{F})= \{0\}$ and by Theorem \ref{orbit}, $E(\mathcal{F})=\{0\}$.
\end{example}
%%%%%%%%%%%%%%%%%%%%%%%%%%%%%%%%%%%%%%%%%%%%%%%%%%%%%%%%%%%%%%%%%%%%%%%%%%%%%%%%%%%%%%%%%%%%%%%%%%%%%%%%%%%%%%%%%%%
For a family $\mathcal{F}$ of entire functions with $F(\mathcal{F})\neq \phi$, the set $$F_{\infty}(\mathcal{F}):=\left\{z\in F(\mathcal{F}):\mbox{ there is a sequence }\left\{f_n\right\}\subset\mathcal{F}\mbox{ such that }f_n(z)\to\infty\right\}$$ is an open as well as closed subset of $F(\mathcal{F})$. Indeed, let  $z_0\in F_{\infty}(\mathcal{F}).$ Then there is a sequence $\left\{f_n\right\}\mbox{ such that }f_n(z_0)\to\infty.$ By normality of $\mathcal{F}\mbox{ at }z_0,$ there is a subsequence $\left\{f_{n_k}\right\}$ of $\left\{f_n\right\}$  which converges uniformly to $\infty$ in some neighborhood $U$ of $z_0$ and hence $U\subset F_{\infty}(\mathcal{F}).$ This proves that $F_{\infty}(\mathcal{F})$ is an open subset of $F(\mathcal{F}).$ Similarly $F_{\infty}(\mathcal{F})$ is closed also.
 
We say that $f\in\partial\left(\mathcal{F}\right)$ if, and only if  there is an open disk $D(z_0,r)\subset F(\mathcal{F})$ and a sequence $\left\{f_n\right\}$ in $\mathcal{F}$ such that $\left\{f_n\right\}$ converges uniformly to $f\mbox{ on }D(z_0,r)$ and $f\notin{\mathcal{F}}.$ By using Vitali's Theorem (\cite{beardon}, p. 56), for a family $\mathcal{F}$ of entire functions, $f\in\partial(\mathcal{F})$ if, and only if there is a sequence $\left\{f_n\right\}\subset\mathcal{F}$ which converges locally uniformly to $f$ on a component of $F(\mathcal{F})$ and $f\notin\mathcal{F}.$
 
 It is observed that $F_{\infty}(\mathcal{F})\neq \phi$ if and only, if $\infty\in \partial\left(\mathcal{F}\right)$. Further,  if $F_{\infty}(\mathcal{F})$ is a non-empty proper subset of $F(\mathcal{F}),$ then $F(\mathcal{F})$ is disconnected. Following example shows that the converse of this statement is not true.
\begin{example}\label{sufficient} Let  $\mathcal{F}_1=\left\{\sin{kz}:k\in\mathbb{N}\right\}.$ Then we show that $J(\mathcal{F}_1)=\mathbb{R}$.

For this, let $z_0\in\mathbb{C}\setminus\mathbb{R}.$ Then chose a disk $D(z_0,r)$ about $z_0$ such that $D(z_0,r)\cap\mathbb{R}=\phi.$ Note that for every $z\in D(z_0,r)\mbox{ and }k\in\mathbb{N},$ $kz\notin\mathbb{R}$ and 
\begin{flalign}
\left|\sin{kz}\right| &=\sqrt{\sin^2kx\cosh^2ky+\cos^2kx\sinh^2ky}\nonumber\\
   &=\sqrt{(1-\cos^2kx)\cosh^2ky+\cos^2kx\sinh^2ky}\nonumber\\
	 &=\sqrt{\cosh^2ky-\cos^2kx\left(\cosh^2ky -\sinh^2ky\right)} \nonumber\\
	 &=\sqrt{\cosh^2ky-\cos^2kx}\nonumber.
\end{flalign}
  Thus $\left|\sin{kz}\right|\to\infty\mbox{ as }k\to\infty$ uniformly on $D(z_0,r).$ Therefore $\mathbb{C}\setminus\mathbb{R}\subset F(\mathcal{F}_1).$ Next, if $z_0\in\mathbb{R}$ then any disk $D(z_0,s)$ about $z_0$ contains a segment of $\mathbb{R}$ which is mapped into $\left[-1,1\right]$ by $\sin{kz}$ for every $k\in\mathbb{N}.$ Whereas for any other point $z\in D(z_0,r)\setminus\mathbb{R},$ $\left|\sin{kz}\right|\to\infty \mbox{ as }k\to\infty.$ So the family $\mathcal{F}_1=\left\{\sin{kz}:k\in\mathbb{N}\right\}$ can not be normal on $z_0\in\mathbb{R}.$ Thus $\mathbb{R}\subset{J(\mathcal{F}_1)}.$ But $\mathbb{C}\setminus\mathbb{R}\subset F(\mathcal{F}_1),$ hence $J(\mathcal{F}_1)=\mathbb{R}.$
 \end{example} 
	
    For $\mathcal{F}_2=\{z^n:n\in\mathbb{N}\}$, $J(\mathcal{F}_2)=\left\{z:|z|=1\right\}$. Let $\mathcal{F}_3=\mathcal{F}_2\cup\mathcal{F}_1.$ Then by Theorem \ref{union} $J(\mathcal{F}_3)=\mathbb{R}\cup\left\{z:|z|=1\right\}.$ Clearly $F(\mathcal{F}_3)$ is disconnected and consists of four components. But 
$F_{\infty}(\mathcal{F}_3)$ is not proper subset of $F(\mathcal{F}_3),$ since it can be easily shown that $F_{\infty}(\mathcal{F}_3)=F(\mathcal{F}_3).$    

%%%%%%%%%%%%%%%%%%%%%%%%%%%%%%%%%%%%%%%%%%%%%%%%%%%%%%%%%%%%%%%%%%%%%%%%%%%%%%%%%%%%%%%%%%%%%%%%%%%%%%%%%%%%%%%%%%%%%%%%%%%%%%%%%%%%%%%%%%%%%%
\section{Proof of main results}
\textbf{Proof of Theorem \ref{union}}:(a) Clearly, $J(\mathcal{F}_1) \cup J(\mathcal{F}_2) \subset J(\mathcal{F}_1 \cup \mathcal{F}_2)$. To show that the other way inclusion, let $z_0\in J(\mathcal{F}_1\cup\mathcal{F}_2).$ Then by Zalcman Lemma \cite{zalcman1},  there is a sequence $\left\{f_n\right\}\subset \mathcal{F}_1\cup\mathcal{F}_2$ , a sequence of positive real numbers $r_n\to 0$ and a sequence $\left\{z_n\right\}:z_n\to z_0 $ as $n\to \infty,$ such that $f_n(z_n+r_nz)$ converges locally uniformly on $\mathbb{C}$ to a non-constant entire function $ f(z).$ There is  a subsequence $\left\{f_{n_k}\right\}$ of $\left\{f_n\right\}$ which lies entirely either in $\mathcal{F}_1\mbox{ or }\mathcal{F}_2$ and $f_{n_k}(z_{n_k}+r_{n_k}z)$ converges locally uniformly on $\mathbb{C}$ to the non-constant entire function $ f(z).$ Hence by the converse to Zalcman Lemma, $z_0\in J(\mathcal{F}_1)\cup J(\mathcal{F}_2).$ Therefore, $J(\mathcal{F}_1\cup\mathcal{F}_2)\subset J(\mathcal{F}_1)\cup J(\mathcal{F}_2).$
 
(b) Suppose that  $\mathbb{C}\setminus U$ contains at least two points. Since $U=\bigcup_{f\in\mathcal{F}}f(N)$, each $f\in\mathcal{F}$ omits at least two distinct values on $N$ and  hence by Montel's theorem, $\mathcal{F}$ is normal in $N,$ which is a contradiction as $z_o\in J(\mathcal{F})\cap N.$  Hence $\mathbb{C}\setminus U$ contains at most one point. $\Box $

\medskip

%%%%%%%%%%%%%%%%%%%%%%%%%%%%%%%%%%%%%%%%%%%%%%%%%%%%%%%%%%%%%%%%%%%%%%%%%%%%%%%%%%%%%%%%%%%%%%%%%%%%%%%%%%%%%%%%%%%%%
\textbf{Proof of Theorem \ref{emptyinterior}}: (a) If $J(\mathcal{F})=D$, then there is nothing to prove. Suppose that $J(\mathcal{F})\neq D$. Assume on the contrary that $Int(J(\mathcal{F}))\neq\phi$. Let $N$ be a neighborhood of some $z\in J(\mathcal{F})$ such that $N\subset J(\mathcal{F})$. Since $J(\mathcal{F})$ is forward invariant, $U=\cup_{f\in\mathcal{F}}f(N)\subset J(\mathcal{F})$. By Theorem \ref{union}, it follows that $\mathbb{C}\setminus J(\mathcal{F})$ contains at most one point. Since $J(\mathcal{F})$ is properly contained in $D$, it follows that $D=\mathbb{C}$ and $U=J(\mathcal{F})$. Since $J({\mathcal{F}})$ is closed in $D=\mathbb{C}$, we have $J(\mathcal{F})=\mathbb{C}=D$, a contradiction. Hence $Int(J(\mathcal{F}))=\phi.$ 
 
To prove $(b),$ suppose  $z_0\in J(\mathcal{F})$ is an isolated point. Then there exists a neighborhood $V$ of $z_0$ such that $V\setminus\left\{z_0\right\}\cap J(\mathcal{F})=\phi.$ Since $f\left(F(\mathcal{F})\right)\subset F(\mathcal{F})$ $\mbox{ for all }f\in\mathcal{F},$  $f\left(V\setminus\left\{z_0\right\}\right)\subset F(\mathcal{F})$ $\mbox{ for all }f\in\mathcal{F}.$ So the family $\mathcal{F}$ omits $J(\mathcal{F})$ on $V\setminus\left\{z_0\right\}$. Therefore, by an extension of  Montel's theorem ( \cite{cara}, p. 203 ), $\mathcal{F}$ is normal in $V,$ which is a contradiction. $\Box$

\medskip

%%%%%%%%%%%%%%%%%%%%%%%%%%%%%%%%%%%%%%%%%%%%%%%%%%%%%%%%%%%%%%%%%%%%%%%%%%%%%%%%%%%%%%%%%%%%%%%%%%%%%%%%%%%%%%%%%%%%%%%%%%%%%%%%%%%%%%%%%%%%%
\textbf{Proof of Theorem \ref{1}}: We use the method of Schwick \cite{schwick} to carry out the proof. Let $f\in\mathcal{F}.$  By an application of the second fundamental theorem of Nevanlinna [\cite{Hayman-1}, p. 44],  the set $A=\left\{w:\Theta(w,f)\geq\frac{1}{2}\right\}$ contains at most two points. Since $J(\mathcal{F})$ contains at least three elements, therefore, for $w_o\in J(\mathcal{F})\setminus A\mbox{, }\Theta(w_o,f)<\frac{1}{2}.$ This implies that the equation $f(z)=w_0$ has infinitely many simple roots $a_1, a_2,...$, say. Now by Zalcman lemma, there is a sequence $f_n\in\mathcal{F},$ a sequence $z_n\to w_0$ and a sequence of positive real numbers $r_n\to 0,$ such that $f_n(z_n+r_nz)\to h(z),$ where $h(z)$ is non-constant entire function. Continuity of $f$ implies that $fof_n(z_n+r_nz)\to foh(z).$
If $h$ is transcendental, then for each $a_n$ except for two values $\Theta (a_n, h)<\frac{1}{2}$ and hence there exists $b \in \mathbb{C}$: $h(b)=a_n$ and $h'(b)\neq 0$. Further, if $h$ is a polynomial, then for each $a_n$ except for one value, $h(z)=a_n$ has simple roots. We pick up one value $a_1$, say, such that there exists $b\in \mathbb{C}$ with $h(b)=a_1$, and $h'(b)\neq 0$ and hence $f(h(b))=w_0$, $f'(h(b)) h'(b)\neq 0$, that is, $(f\circ h)^{'}(b)\neq 0.$
   	  Next, $fof_n(z_n+r_nz)-(z_n+r_nz)\to foh(z)-w_0.$ Since $foh-w_0$ has zero at $z=b$ and $foh-w_0$ is not constant, by Hurwitz theorem,  $fof_n(z_n+r_nz)-(z_n+r_nz)$ has zeros at $c_n \mbox{ with }c_n\to b$ for all sufficiently large $n$. Thus $w_n=z_n+r_nc_n$ is a fixed point of $fof_n.$ Since, for large $n, r_n\left(fof_n\right)'(z_n+r_nc_n)=\left(fof_n(z_n+r_nz)\right)'(c_n)\to (foh)'(b)\neq 0$ so that $\left|\left(fof_n\right)'(z_n+r_nc_n)\right|>1.$ $\Box$

\medskip
	
%%%%%%%%%%%%%%%%%%%%%%%%%%%%%%%%%%%%%%%%%%%%%%%%%%%%%%%%%%%%%%%%%%%%%%%%%%%%%%%%%%%%%%%%%%%%%%%%%%%%%%%%%%%%%%%%%%%%%
 The proof of Theorem \ref{poly} is on the similar lines as that of Theorem \ref{1}.

\medskip

%%%%%%%%%%%%%%%%%%%%%%%%%%%%%%%%%%%%%%%%%%%%%%%%%%%%%%%%%%%%%%%%%%%%%%%%%%%%%%%%%%%%%%%%%%%%%%%%%%%%%%%%%%%%%%%%%%%%%
\textbf{Proof of Theorem \ref{orbit}}:(a) If $J(\mathcal{F})=\phi $, then there is nothing to prove. Suppose that $J(\mathcal{F})\neq\phi.$ Assume the contrary that there is $z_0\in J({\mathcal{F}})$ such that $z_0\notin \overline{\mathcal{O}^{-}_{\mathcal{F}}(z_1)}$ for some $z_1 \in \mathbb{C} \setminus E(\mathcal{F})$. That is, there is a neighborhood $N$ of $z_0$ such that $N\cap \mathcal{O}^{-}_{\mathcal{F}}(z_1)=\phi$. We choose a neighborhood $N_1 \subset N$ of $z_0$ such that $\left(N_1 \setminus{\{z_0\}}\right)\cap\mathcal{O}^{-}_{\mathcal{F}}(z_2)=\phi$ for some $z_2\in E(\mathcal{F})$ since $\mathcal{O}^{-}_{\mathcal{F}}(z_2)$ is a finite set. Then $\cup_{f\in\mathcal{F}}f(N_{1}\setminus\{z_0\})$ omits  $z_1$ and $ z_2$. Therefore, by an extension of  Montel's theorem ( \cite{cara}, p.203 ), $\mathcal{F}$ is normal in $N_1$, which is a contradiction as $z_0\in J(\mathcal{F}) \cap N_1$.
 
(b) Suppose that $E(\mathcal{F})\geq 2$ and let $z_1,z_2\in E(\mathcal{F})$. Let $z_0\in J(\mathcal{F})$. Since $\mathcal{O}^{-}_{\mathcal{F}}(z_1)\cup \mathcal{O}^{-}_{\mathcal{F}}(z_2)$ is a finite set, we choose a neighborhood $N$ of $z_0$ such that $N\setminus{\{z_0\}}\cap\left(\mathcal{O}^{-}_{\mathcal{F}}(z_1)\cup \mathcal{O}^{-}_{\mathcal{F}}(z_2)\right)=\phi$ and hence by an extension of Montel's theorem, $\mathcal{F}$ is normal in $N$, which is a contradiction as $z_0\in J(\mathcal{F})\cap N.$ $\Box$

\medskip

\bibliographystyle{amsplain}

%%%%%%%%%%%%%%%%%%%%%%%%%%%%%%%%%%%%%%%%%%%%%%%%%%%%%%%%%%%%%%%%%%%%%%%%%%%%%%%%%%%%%%%%%%%%%%%%%%%%%

%%%%%%%%%%%%%%%%%%%%%%%%%%%%%%%%%%%%%%%%%%%%%%%%%%%%%%%%%%%%%%%%%%%%%%%%%%%%%%%%%%%%%%

\end{document}